%% file: transient_thermal_problems_on_the_half-space.tex
\include{macros}

\documentclass[a4paper,12pt]{article}

\usepackage{amsmath}
\usepackage{amssymb}
\usepackage{verbatim}
\usepackage{graphicx,color}
\usepackage{psfrag}
\author{William J.~Parnell$^{*}$, Vu-Hieu Nguyen$^{\dagger}$, Raphael Assier$^*$,\\Salah Naili$^{\dagger}$ and I.\ David Abrahams$^{*}$\\
\footnotesize{$*$ School of Mathematics, University of Manchester, Oxford Road, Manchester, M13 9PL, UK}\\
\footnotesize{$\dagger$ Universit\'{e} Paris-Est, Laboratoire Modelisation et Simulation Multi Echelle,}\\ \footnotesize{MSME UMR 8208, CNRS, 61 avenue du G\'{e}n\'{e}ral de Gaulle, 94010 Cr\'{e}teil Cedex, France}}
\title{Transient thermal mixed boundary value problems in the half-space}

\begin{document}

\maketitle

\numberwithin{equation}{section}
\begin{abstract}
The Wiener-Hopf and Cagniard-de Hoop techniques are employed in order to solve a range of transient thermal mixed boundary value problems on the half-space. The thermal field is determined via a rapidly convergent integral, which can be evaluated straightforwardly and quickly on a desktop PC.

\end{abstract}

\subsection*{Keywords}

Heat conduction; Transient; Mixed boundary conditions; Wiener-Hopf; Cagniard-de Hoop

\section{Introduction} \label{intro}

Traditionally, great interest has been shown in determining the disturbances that are generated when loads are applied on the surface of a half-space. Lamb \cite{Lam-04} obtained the exact solution when an impulsive, concentrated load is applied along a line of the free surface of an isotropic linear elastic medium. de Hoop reappraised this problem \cite{DeH-60}, modifying the method originally devised by Cagniard \cite{Cag-39}, \cite{Cag-62}, leading to the now well-known \textit{Cagniard-de Hoop} (CdH) technique. This method has been used widely since, allowing exact solutions to be obtained for a wide range of transient elasticity problems. The method can also be useful in order to render solutions into integral forms that are rapidly convergent when calculated numerically.

Transient \textit{thermoelastic} half-space problems were considered by Danilovskaya \cite{Dan-50}, Boley and Tolins \cite{Bol-62} and Achenbach \cite{Ach-63} but in these problems the forcing was such that the CdH technique was not required. The extension of these problems to inhomogeneous media was considered by Baczynski \cite{Bac-03} and Parnell \cite{Par-06}. A purely thermal, transient problem that employed the CdH method was solved in \cite{She-01}. The thermoelastic Lamb problem was studied by Nayfeh and Nemat-Nasser \cite{Nay-72} who used generalized thermoelasticity in order to retain a finite thermal wave speed, employing the CdH technique to determine the solution.

All of the above problems are of fundamental importance in an array of applications where a number of alternative boundary conditions on the surface can arise. What appears to be rather lacking in the literature however are studies of transient problems with \textit{mixed} boundary conditions, where in the context of the thermal problems the condition takes the form, e.g.
\begin{align}
\cT(0,y,t) &= f(y,t) & \textnormal{for $y>0$}, &&
\deriv{\cT}{x}(0,y,t) &= g(y,t) &&  \textnormal{for $y<0$}, \label{mixed1}
\end{align}
where $\cT(x,y,t)$ is the temperature field, $f(y,t)$ and $g(y,t)$ are two specified functions, $t$ is time and with reference to Fig.\ \ref{fig:setdomaindescription} $x$ and $y$ are Cartesian coordinates. The half-space resides in $x\geq 0$ and $y$ runs parallel to the surface, which is defined by $x=0$.

Generally such problems lead to the propagation of a thermal disturbance into the half-space. Indeed, such thermal front problems are of importance in a number of applications including defect sizing \cite{Alm-94}, transient thermography \cite{Sha-13}, solar cell manufacturing \cite{Pil-02} and thermal insulation \cite{Ree-00}. Caflisch and Keller \cite{Caf-81},  Levine \cite{Lev-82} and Satapathy and Sahoo \cite{Sat-02} studied front propagation in the thermal context with mixed boundary conditions but in the context of steady problems with applications in quenching. Kozlov et al.\ \cite{Koz-01} considered a transient half-space problem with mixed boundary conditions and made progress by using cylindrical coordinates due to the special form of the boundary condition chosen.

Mixed boundary conditions are generally difficult to handle even in steady problems and analytical or semi-analytical solutions are frequently only possible by the application of the Wiener-Hopf method \cite{Nob-88}. This method exploits the analyticity properties of functions in order to yield an explicit or approximate solution in the Fourier transform domain. Contour integration then yields the solution in the physical domain.


Here we shall consider a rather general mixed boundary value problem in the context of thermal front propagation and determine solutions using the Wiener-Hopf method and Cagniard-de Hoop technique. This problem of mixed boundary conditions of the form \eqref{mixed1} is of particular interest in analyzing the field close to the location of the change in boundary condition type, i.e.\ $x=y=0$ in \eqref{mixed1}.

We obtain a solution in single integral form by using a deformation of the Laplace contour in a similar manner to the Cagniard-de Hoop method. Although it appears that we cannot obtain an explicit solution, the solution determined can be evaluated rapidly on a desktop PC and therefore it is of great utility due to its general form and its ability to circumvent a direct numerical simulation of the problem. Although similar problems, involving a discontinuous temperature boundary condition have been considered in the building insulation literature, see e.g.\ Claesson and Hegentoft \cite{Cla-91} and Hegentoft and Claesson \cite{Heg-91} to the authors' knowledge it does not appear that the solution we provide has been written down anywhere in the literature before now.

In this paper we shall first set out the problem description in section \ref{problemdescription} before determining the solution in the transform domain in section \ref{WH}. In section \ref{CDH} we describe how we deform the Laplace contour onto a steepest descent path, in the manner of the Cagniard-de Hoop technique in order to obtain a solution in terms of a single integral along the deformed contour path, with an integrand that decays exponentially. In section \ref{sec:particularexamples} we illustrate the efficacy of the scheme by determining the solution for a number of different boundary conditions, with validation provided by finite element solutions.

\section{Problem description} \label{problemdescription}

Assume that the problem under consideration is two-dimensional, being independent of $z$ and define the two dimensional half-space domain $\mathcal{D}=\{(x,y):0\leq x <\infty, -\infty < y < \infty\}$. We seek solutions to the anisotropic heat equation:
\begin{align}
\frac{k}{\rho c_V}\left(\derivtwo{\cT}{x}+\ell\derivtwo{\cT}{y}\right) &= \deriv{\cT}{t}
\end{align}
where $k$ and $k\ell$ are the thermal conductivities ($\ell>0$) in the $x$ and $y$ directions respectively, $c_V$ is the specific heat at constant volume, $\rho$ is the mass density, $t$ is time and $\cT=\cT(x,y,t)$ is the temperature field. We can combine the constants as $\ka=k/(\rho c_V)$, the thermal diffusivity.

\begin{figure}[h!]
\begin{center}
\includegraphics[scale=0.6]{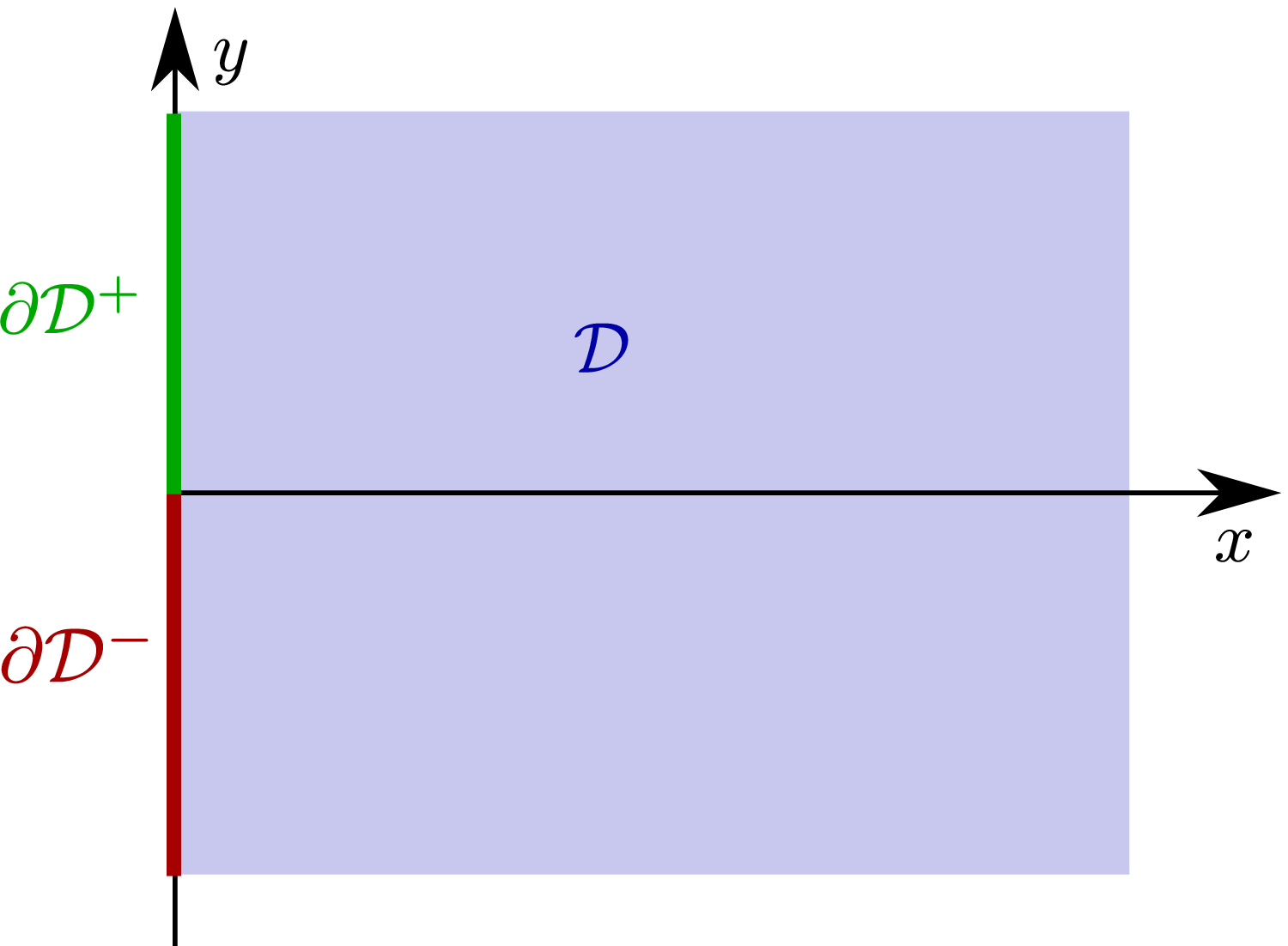}
\end{center}
\caption{Domain $\mathcal{D}$ of the problem and its boundaries $\partial \mathcal{D}^+$ and $\partial \mathcal{D}^-$} on which Dirichlet and Neumann boundary conditions are imposed respectively.
\label{fig:setdomaindescription}
\end{figure}

It is convenient to non-dimensionalise the governing equation, using coordinates with a ``hat'' and scale the $y$ coordinate to remove the anisotropy coefficient. Write $(\hat{x},\hat{y},\hat{T},\hat{t})=(x/x^*,y/y^*,(\cT-\cT^*)/\cT^*,t/t^*)$, where
\begin{align}
x^*=1[\text{m}], \quad y^*=\frac{1}{\sqrt{\ell}} [\text{m}], \quad t^*=\frac{(x^*)^2}{\ka} [\text{s}],
\end{align}
and $\cT^*$ is the reference temperature in Kelvin. Upon doing so and ``dropping hats'' we find
\begin{align}
\nabla^2 \cT &= \deriv{\cT}{t}    \label{heat}
\end{align}
where $\nabla^2 = \pa^2/\pa x^2 + \pa^2/\pa y^2$. We wish to solve \eqref{heat} on the (scaled) domain $\mathcal{D}$ with boundary $\partial \mathcal{D}=\partial \mathcal{D}^-\cup\partial \mathcal{D}^+$ as illustrated in Fig.\ \ref{fig:setdomaindescription}. We consider homogeneous initial conditions of the form
\begin{align}
\cT(x\geq 0,y,t=0) &= 0 \label{2.4}
\end{align}
and boundary conditions of the form \eqref{mixed1} but simplify by removing the $y$ dependence, i.e.\
\begin{align}
\cT\Big|_{\partial \mathcal{D}^+}=\cT(x=0,y>0,t>0)= T_0 f_0(t), \\
\deriv{\cT}{x}\Big|_{\partial \mathcal{D}^-}= \deriv{\cT}{x}(x=0,y<0,t>0) = T_0'g_0(t), \label{2.6}
\end{align}
where $T_0$ and $T_0'$ are real constants and $f_0(t)$ and $g_0(t)$ are piecewise continuous functions of time.

We therefore have a mixed boundary value problem, which in general are not straightforward to solve even in the steady context so that the time dependence adds an additional element of complexity. Furthermore we allow for the fact that we could have a step change at $t=0$ on $x=0$, leading to a propagating discontinuity front in the half-space.

In order to determine $\cT$ it is convenient to introduce an alternative problem (for convergence issues as will be shown), giving rise to a different temperature distribution $T$, depending on a small parameter $\epsilon$ and where $T$ converges to $\cT$ as $\epsilon\rightarrow 0$. This problem is described as follows
\begin{align}
\nabla^2 T &= \deriv{T}{t}\label{eq:newheateq}\\
T(x\geq0,y,t=0)&=0\\
T\Big|_{\partial \mathcal{D}^+}=T(x=0, y>0,t>0)&=T_0f_0(t)e^{-\epsilon y}\\
 \deriv{T}{x}\Big|_{\partial \mathcal{D}^-}=\deriv{T}{x}(x=0,y<0,t>0)&=T_0'g_0(t)e^{\epsilon y}. \label{eq:newderbc}
\end{align}
As is easily seen, we recover the solution to our original problem by taking the limit as $\epsilon$ tends to zero:
\begin{align}
\cT(x,y,t)=\lim_{\epsilon \rightarrow 0} T(x,y,t).
\end{align}

\section{Solution in the transform domains via the Wiener-Hopf technique} \label{WH}

Define the Laplace transform in time for any function $\phi(x,y,t)$ by
\begin{align}
\mathcal{L}(\phi(x,y,t)) = \tilde{\phi}(x,y,s) &= \int_0^{\infty}\phi(x,y,t)e^{-st}{dt}
\end{align}

and hence applying this to the governing scaled equations (\ref{eq:newheateq})-(\ref{eq:newderbc}) we have
\begin{align}
\nabla^2 \tilde{T} &= s\tilde{T} \label{heatLT}
\end{align}
and boundary conditions become
\begin{align}
\tilde{T}(x=0,y>0,s) &= \tilde{f}_0(s)T_0e^{-\epsilon y} \\
\deriv{\tilde{T}}{x}(x=0,y<0,s) &= \tilde{g}_0(s)T_0'e^{\epsilon y}.
\end{align}
Although $s\in\mathbb{C}$ the set of complex numbers, for the sake of the analysis to follow we can assume it to be real and positive. This allows us to scale the $(x,y)$ variables to simplify the governing equations. The derivation goes through retaining explicit dependence on $s$ in the governing equation, but the algebra becomes rather heavy and tedious and does not render any greater understanding of the problem; both approaches lead to the same result. As such we will rescale $x$ and $y$ in order to eliminate $s$ from the governing equation. Thus define
\begin{align}
x_0 &= x\sqrt{s}, & y_0 &= y\sqrt{s} \label{x0y0}
\end{align}
and therefore we obtain
\begin{align}
\nabla^2_0 \tilde{T} &= \tilde{T} \label{heatLT2}
\end{align}
where $\nabla_0^2 = \pa^2/\pa x_0^2 + \pa^2/\pa y_0^2$. The boundary conditions become
\begin{align}
\tilde{T}(x_0=0,y_0>0,s) &= \tilde{f}_0(s)T_0e^{-\epsilon y} \\
\deriv{\tilde{T}}{x_0}(x_0=0,y_0<0,s) &= \frac{\tilde{g}_0(s)}{\sqrt{s}}T_0'e^{\epsilon y}.
\end{align}

\begin{figure}[h!]
\begin{center}
\includegraphics[scale=0.6]{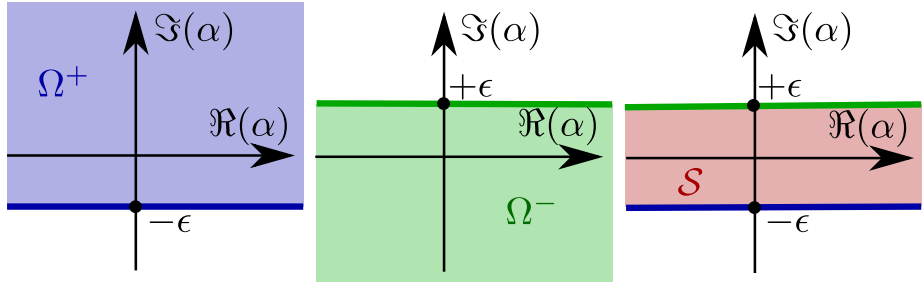}
\end{center}
\caption{Illustrating the domains $\Omega^+$, $\Omega^-$ and $\mathcal{S}$ on which the functions $\Theta^+, \Theta^-$ and $\Theta$ are analytic, respectively.}
\label{fig:setdescription}
\end{figure}

Next define the Fourier transform in $y_0$ as
\begin{align}
\mathcal{F}(\tilde{T}(x_0,y_0,s)) = \Theta(x_0,\al,s) &= \int_{-\infty}^{\infty}\tilde{T}(x_0,y_0,s)e^{i\al y_0}{dy_0}
\end{align}
and define $\Theta^+$ and $\Theta^-$ as
\begin{align}
\Theta^-(x_0,\al,s)= \int_{-\infty}^{0}\tilde{T}(x_0,y_0,s)e^{i\al y_0}{dy_0} \quad \text{and} \quad \Theta^+(x_0,\al,s)=\int_{0}^{\infty}\tilde{T}(x_0,y_0,s)e^{i\al y_0}{dy_0},
\end{align}
so that $\Theta=\Theta^-+\Theta^+$. Applying the Fourier transform to the governing (Laplace transformed) equation \eqref{heatLT2} we find that
\begin{align}
\Theta'' &= (\al^2+1)\Theta \label{heattrans}
\end{align}
and to the boundary conditions, we find that
\begin{align}
\Theta^+(x_0=0,\al,s) &= \frac{i T_0}{(\al+i\eps)}\tilde{f}_0(s) \label{BC1}\\
\deriv{\Theta^-}{x_0}(x_0=0,\al,s) &= -\frac{i T_0'}{(\al-i\eps)}\frac{\tilde{g}_0(s)}{\sqrt{s}} \label{BC2}.
\end{align}
With reference to Fig.\ \ref{fig:setdescription}, we note that $\Theta^+$ is analytic on $\Omega^+=\{\al \in \mathbb{C}, \Im(\al)>-\epsilon\}$, $\Theta^-$ is analytic on $\Omega^-=\{\al \in \mathbb{C}, \Im(\al)<\epsilon\}$ and $\Theta^+$, $\Theta^-$ and $\Theta$ are analytic on the strip $\mathcal{S}=\Omega^- \cap \Omega^+$. The superscript $+$ and $-$ notation thus indicates analyticity in the domains $\Omega^+$ and $\Omega^-$ respectively.

The solution of \eqref{heattrans} is
\begin{align}
\Theta(x_0,\al,s) &= A_1(\al,s) \exp(-(\al^2+1)^{1/2}x_0) \label{Thetasol}
\end{align}
where the branch of the square root function in the exponent is chosen so that its real part satisfies $\Re(\al^2+1)^{1/2}>0$, ensuring that the solution decays as $x_0\rightarrow \infty$.  For conciseness let us introduce $\Phi^+$ and $\Psi^-$, analytic on $\Omega^+$ and $\Omega^-$ respectively, as
\begin{align}
\Phi^+(\al,s)=\deriv{\Theta^+}{x}\Big|_{x=0} \quad \text{and} \quad \Psi^-(\al,s)=\Theta^{-}\Big|_{x=0}. \label{concise}
\end{align}
Imposing the boundary conditions \eqref{BC1}-\eqref{BC2}, employing \eqref{concise} and eliminating $A_1$ between the two resulting equations, we arrive at
\begin{align}
-\Phi^+ &= K(\al)\Psi^- + \frac{K(\al)iT_0}{(\al+i\eps)}\tilde{f}_0(s)-\frac{iT_0'}{(\al-i\eps)}\frac{\tilde{g}_0(s)}{\sqrt{s}}. \label{elimA1}
\end{align}
The \textit{kernel} here $K(\al)=(\al^2+1)^{1/2}$ is easily factorized as $K(\al)=K^-/K^+$ where
\begin{align}
K^+(\al) = (\al+i)^{-1/2} \quad \text{and} \quad  K^-(\al) = (\al-i)^{1/2}. \label{Kminus}
\end{align}

Multiplying \eqref{elimA1} by $K^+$  we obtain
\begin{align}
-K^+(\al)\Phi^+ &= K^-(\al)\Psi^- + S(\al,s) \label{elimA12}
\end{align}
where
\begin{align}
S(\al,s) &= \frac{K^-(\al)iT_0}{(\al+i\eps)}\tilde{f}_0(s)-\frac{K^+(\al)iT_0'}{(\al-i\eps)}\frac{\tilde{g}_0(s)}{\sqrt{s}}  \label{S} \\
 &= S^-(\al,s) + S^+(\al,s),
\end{align}
where we have indicated that we wish to determine a sum factorization of the function $S$.  One can employ the pole removal method \cite{Vei-07} in order to show quite straightforwardly that
\begin{align}
S^+(\al,s)&=iT_0f_0(s)L_1^++iT_0'\frac{g_0(s)}{\sqrt{s}}L_2^+, \\
S^-(\al,s)&=iT_0f_0(s)L_1^-+iT_0'\frac{g_0(s)}{\sqrt{s}}L_2^-,\label{Sminus}
\end{align}
where
\begin{align}
L_1^- &= \frac{(\al-i)^{1/2}-c_-}{\al+i\eps}, &  \quad L_1^+=\frac{c_-}{\al+i\eps}, \label{eq:L1}, \\
L_2^+ &= -\frac{(\al+i)^{-1/2}-c_+}{\al-i\eps}, & L_2^-=-\frac{c_+}{\al-i\eps} \label{eq:L2}
\end{align}
and
\begin{align}
c_- &= [-i(\eps+1)]^{1/2}, & c_+ &= [i(\eps+1)]^{-1/2}.
\end{align}

Referring to \eqref{elimA1} we can therefore define a function $E(\al)$ such that
\begin{align}
E(\al)=\left\{ \begin{array}{ccc}
-(K^+\Phi^++S^+)&\text{on}&\Omega^+\backslash \mathcal{S} \\
-(K^+\Phi^++S^+)=(K^-\Psi^-+S^-)&\text{on}&\mathcal{S}\\
(K^-\Psi^-+S^-)&\text{on}&\Omega^-\backslash \mathcal{S},
\end{array} \right. \label{EE}
\end{align}
and therefore is analytic on the whole $\alpha$-complex plane. It can be shown that as $|\al|\rightarrow\infty$, $E(\alpha)=O(\al^{-1/2})$ when $\al\in\Omega^+$ and $E(\al)=O(\al^{1/2})$ when $\al\in\Omega^-$. As such $E(\al)=o(\al)$ as $|\al|\rightarrow\infty$. This together with the analyticity of $E(\al)$ and the extended Liouville theorem (see for example \cite{Nob-88}), implies that $E(\al)$ is constant. However, we also know that $E(\al)\rightarrow 0$ as $|\al|\rightarrow\infty$ and $\al\in\Omega^+$. We can then conclude that $E(\al)=0$ everywhere. Given this, we therefore have from \eqref{EE}
\begin{align}
\Phi^+ = -\frac{S^+}{K^+}\quad \text{and} \quad \Psi^- &= -\frac{S^-}{K^-}.
\end{align}
From the original expressions for $A_1$ determined from the boundary conditions we can show that
\begin{align}
A_1(\al,s) &= \frac{iT_0}{(\al+i\eps)}\tilde{f}_0(s) +\Psi^- \\ 
&= \frac{iT_0 c_- \tilde{f}_0(s)}{(\al+i\eps)(\al-i)^{1/2}} + \frac{iT_0' c_+ \tilde{g}_0(s)}{(\al-i\eps)(\al-i)^{1/2}\sqrt{s}}
\end{align}
where  \eqref{Kminus}, \eqref{eq:L1}, \eqref{eq:L2} and \eqref{Sminus} have all been used. Referring to \eqref{Thetasol}, the solution in transform space is therefore
\begin{align}
\Theta(x_0,\al,s) 
 &= \left(\frac{iT_0 f \tilde{f}_0(s)}{(\al+i\eps)(\al-i)^{1/2}} + \frac{iT_0' g \tilde{g}_0(s)}{(\al-i\eps)(\al-i)^{1/2}\sqrt{s}}\right)\exp(-(\al^2+1)^{1/2}x_0).
\end{align}


Formally inverting the Fourier transform and using \eqref{x0y0} gives the Laplace transformed solution as
\begin{align}
\tilde{T}(x,y,\al) &= \frac{iT_0\tilde{f}_0(s)}{2\pi}I_1(x,y,s)+\frac{iT_0'\tilde{g}_0(s)}{2\pi\sqrt{s}}I_2(x,y,s), \label{eq:defI1I2raph}
\end{align}
where
\begin{align}
I_1=c_-\int_{-\infty}^{\infty}\frac{e^{-[(\al^2+1)^{1/2}x+i\al y]\sqrt{s}}}{(\al+i\eps)(\al-i)^{1/2}} \, \mathrm{d}\al \quad \text{and} \quad I_2=c_+\int_{-\infty}^{\infty}\frac{e^{-[(\al^2+1)^{1/2}x+i\al y]\sqrt{s}}}{(\al-i\eps)(\al-i)^{1/2}} \, \mathrm{d}\al \label{I1I2}
\end{align}

\section{Semi-analytical inversion via a Cagniard-de Hoop approach} \label{CDH}

\begin{figure}[h!]
\begin{center}
\includegraphics[scale=0.6]{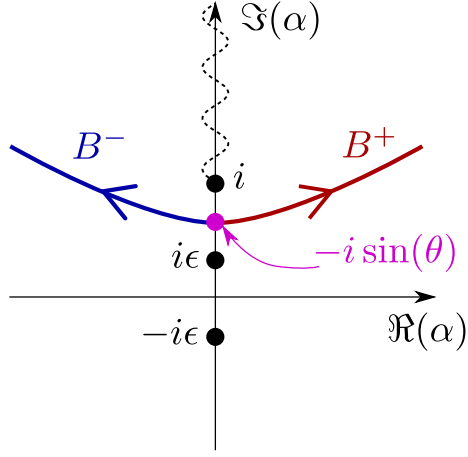}
\end{center}
\caption{Illustrating the location of the paths $B^\pm$ when $\theta\in[-\pi/2,0)$.  The poles at $\al=\pm i\eps$ and branch point at $\al=i$ are also depicted where the branch cut passes to infinity vertically along the imaginary axis.}
\label{fig:alphaplane1}
\end{figure}

\begin{figure}[h!]
\begin{center}
\includegraphics[scale=0.6]{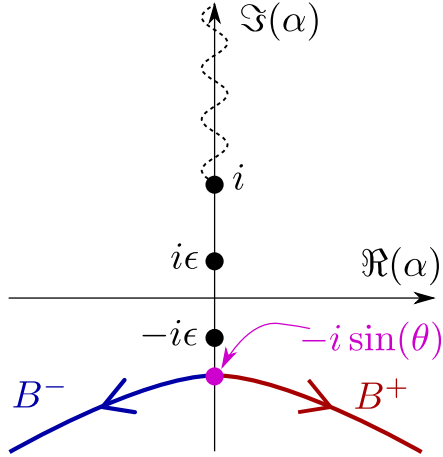}
\end{center}
\caption{Illustrating the location of the paths $B^\pm$ when $\theta\in(0,\pi/2]$. The poles at $\al=\pm i\eps$ and branch point at $\al=i$ are also depicted where the branch cut passes to infinity vertically along the imaginary axis.}
\label{fig:Blowerraph}
\end{figure}

Motivated by the Cagniard-de Hoop technique, let us introduce polar coordinates $r$ and $\theta$ related to $x$ and $y$ in the usual manner, i.e.\ $x=r\cos\theta, y=r\sin\theta$ where $\theta\in[-\pi/2,\pi/2]$ and introduce the parameter $\beta$ via the expression
\begin{align}
\beta r &= (\al^2+1)^{1/2} x + i\al y
\end{align}
so that
\begin{align}
\beta &= (\al^2+1)^{1/2} \cos\theta + i\al \sin\theta.
\end{align}
Inverting for $\alpha$ we therefore determine the two paths $B_+$ and $B_-$ in the right and left halves of the $\alpha$-plane
\begin{align}
\al_{\pm} &= -i\be\sin\theta \pm \sqrt{\be^2-1}\cos\theta.
\end{align}
In the $\alpha$-plane, with $\be\in[1,\infty)$ these paths start at $\al=-i\sin\theta$ and move off to infinity either in the upper half-plane ($\theta\in(-\pi/2,0)$, see Fig.\ \ref{fig:alphaplane1}) or the lower half-plane ($\theta\in(0,\pi/2)$, see Fig.\ \ref{fig:Blowerraph}).

Since $B^{\pm}$ are steepest descent paths for the integrals, the idea is to deform the integrals \eqref{I1I2} from the real line onto these to aid convergence. In classical Cagniard de Hoop problems this frequently permits one to write the $\al$ integral in the form of a Laplace transform of a function that is independent of $s$ and thus we can determine the inverse transform immediately, thus rendering explicit solutions. Here we are not so fortunate, the function will not be independent of $s$ but nevertheless we are able to make significant progress due to the fact that the inverse Laplace transform integral can be determined analytically in many important cases as we shall see shortly. This leaves the solution in a single integral form that is rapidly convergent.

At this point note that
\begin{align}
\lim_{\eps\rightarrow 0} c_- = (e^{-i\pi/2})^{1/2} = e^{-i\pi/4} = c \quad \text{and} \quad
\lim_{\eps\rightarrow 0} c_+ = (e^{i\pi/2})^{-1/2} = e^{-i\pi/4} = c
\end{align}
and let us consider the case of negative and positive $\theta$ separately.

\subsection{The case of $\theta\in[-\pi/2,0)$} \label{sec:negtheta}

\subsubsection{Evaluation of $I_1$}

Referring to \eqref{I1I2} and Fig.\ \ref{fig:alphaplane1}, we see that the integrand of $I_1$ has a pole at $\al=-i\eps$ and a branch point at $\al=i$. Deforming the contour from the real axis to the $B^\pm$ contour in the $\al$-plane we do not cross any singularities and hence we find that
\begin{align}
I_1 
&=c_-\int_{B^+}\frac{e^{-[(\al^2+1)^{1/2}x+i\al y]\sqrt{s}}}{(\al+i\eps)(\al-i)^{1/2}}\, \mathrm{d}\al-c_-\int_{B^-}\frac{e^{-[(\al^2+1)^{1/2}x+i\al y]\sqrt{s}}}{(\al+i\eps)(\al-i)^{1/2}}\, \mathrm{d}\al \nonumber\\
&=I_1^+-I_1^- \label{I1s}
\end{align}
On $B^\pm$, we have $\al=\al_\pm=-i\be\sin\theta \pm \sqrt{\be^2-1}\cos\theta$ and so we have
\begin{align}
\mathrm{d}\al=(-i\sin\theta\pm\be\cos\theta(\be^2-1)^{-1/2})\mathrm{d}\be.
\end{align}
Noting that $\be\in[1,\infty)$ is a parametrization of the paths $B^{\pm}$ we can rewrite $I_1^\pm$ as follows, taking the limit as $\eps\rightarrow 0$,
\begin{align}
\lim_{\eps \rightarrow 0} I_1^\pm=c\int_1^\infty  \frac{(-i\sin\theta\pm\be\cos\theta(\be^2-1)^{-1/2})}{\al_\pm(\al_\pm-i)^{1/2}} e^{-\be r \sqrt{s}} \,\mathrm{d}\be
\end{align}
Therefore from \eqref{I1s} we determine the form
\begin{align}
\lim_{\eps \rightarrow 0} I_1 
&= c\int_1^\infty \mathcal{F}(\be,\theta) e^{-\be r \sqrt{s}} \,\mathrm{d}\be, \label{eq:limL1raph}
\end{align}
where
\begin{align}
\mathcal{F}(\be,\theta)=&-i\sin\theta\left(\frac{1}{\al_+(\al_+-i)^{1/2}}-\frac{1}{\al_-(\al_--i)^{1/2}} \right)\nonumber \\
&+\frac{\be}{(\be^2-1)^{1/2}}\cos\theta\left( \frac{1}{\al_+(\al_+-i)^{1/2}}+\frac{1}{\al_-(\al_--i)^{1/2}}\right). \label{eq:calFraph}
\end{align}

As an aside, we note by deforming into the lower half-plane that
\begin{align}
\int_1^{\infty}\mathcal{F}(\be,\theta)\,\mathrm{d}\be &=\lim_{\epsilon \rightarrow 0} \int_{-\infty}^{\infty} \frac{1}{(\al+i \epsilon)(\al-i)^{1/2}}\,\mathrm{d}\al\\
 &= - \frac{2i\pi}{c}. \label{eq:easynegraph}
\end{align}

\subsubsection{Evaluation of $I_2$}

Referring to \eqref{I1I2} and Fig.\ \ref{fig:alphaplane1}, we see that the integrand of $I_2$ has a simple pole at $\al=i\eps$ and a branch point at $\al=i$. Deforming the contour from the real axis to the $B^\pm$ contour in the $\al$-plane we cross the simple pole and pick up its residue. Accounting for this contribution we find that
\begin{align}
I_2 
&= c_+\int_{B^+} \frac{e^{-[(\al^2+1)^{1/2}x+i\al y]\sqrt{s}}}{(\al-i\eps)(\al-i)^{1/2}} \, \mathrm{d}\al-c_+\int_{B^-} \frac{e^{-[(\al^2+1)^{1/2}x+i\al y]\sqrt{s}}}{(\al-i\eps)(\al-i)^{1/2}} \, \mathrm{d}\al+R_2\nonumber\\
&=I_2^+-I_2^-+R_2 \label{I2s}
\end{align}
where
\begin{align}
R_2 
&=2\pi i c_+\frac{e^{-(((i\eps)^2+1)^{1/2} \cos\theta + i(i\eps) \sin\theta)r\sqrt{s}}}{(i\eps-i)^{1/2}}
\end{align}
and in particular,
\begin{align}
\lim_{\eps \rightarrow 0} R_2 
&=2i\pi e^{-x\sqrt{s}}. \label{R2s}
\end{align}
Noting that $\be\in[1,\infty)$ is a parametrization of the paths $B^{\pm}$ we can rewrite $I_2^\pm$ as follows, taking the limit as $\eps\rightarrow 0$,
\begin{align}
\lim_{\eps \rightarrow 0} I_2^\pm= c\int_1^\infty  \frac{(-i\sin\theta\pm\be\cos\theta(\be^2-1)^{-1/2})}{\al_\pm(\al_\pm-i)^{1/2}} e^{-\be r \sqrt{s}} \,\mathrm{d}\be.
\end{align}
Therefore from \eqref{I2s} and \eqref{R2s} we determine the form
\begin{align}
\lim_{\eps \rightarrow 0} I_2 
&= c\int_1^\infty \mathcal{F}(\be,\theta) e^{-\be r \sqrt{s}} \,\mathrm{d}\be+2\pi i e^{-x\sqrt{s}}. \label{eq:limL2raph}
\end{align}

\subsubsection{An expression for $\mathcal{T}(r,\theta,t)$}

We show in Appendix \ref{appendix} that
\begin{align}
c\mathcal{F}(\be,\theta)=\frac{\sqrt{2}}{i}\mathcal{G}(\be,\theta), \label{eq:calGraph}
\end{align}
where $\mathcal{G}$ is a real-valued function. As such, using this together with \eqref{eq:defI1I2raph}, \eqref{eq:limL1raph}, \eqref{eq:limL2raph} and \eqref{eq:calGraph} we find the following expression for $\tilde{\mathcal{T}}=\lim_{\eps \rightarrow 0}\tilde{T}$
\begin{align}
\tilde{\mathcal{T}}(r,\theta,s)&=\frac{iT_0\tilde{f}_0(s)}{2\pi}\lim_{\eps \rightarrow 0}I_1(x,y,s)+\frac{iT_0'\tilde{g}_0(s)}{2\pi\sqrt{s}}\lim_{\eps \rightarrow 0}I_2(x,y,s)\nonumber\\
&= \frac{1}{\sqrt{2}\pi}\left(T_0\tilde{f}_0(s) + T_0'\frac{\tilde{g}(s)}{\sqrt{s}}\right) \int_1^\infty \mathcal{G}(\be,\theta) e^{-\be r \sqrt{s}} \,\mathrm{d}\be -\frac{T_0'\tilde{g}_0(s)}{\sqrt{s}} e^{-x\sqrt{s}}.
\end{align}
Finally we recover $\mathcal{T}$ by taking the inverse Laplace Transform,
\begin{align}
\mathcal{T}(r,\theta,t) 
&=\frac{T_0}{\sqrt{2}\pi}\int_1^\infty \mathcal{G}(\be,\theta) \mathcal{T}_1(r,\be,t)\, \mathrm{d}\be+\frac{T_0'}{\sqrt{2}\pi}\int_1^\infty\mathcal{G}(\be,\theta)\mathcal{T}_2(r,\be,t)\, \mathrm{d}\be \nonumber \\
& \hspace{4cm}-T_0'\mathcal{T}_2(r,\cos\theta,t),
\end{align}
where
\begin{align}
\mathcal{T}_1(r,\be,t) &= \frac{1}{2i\pi}\int_{\si-i\infty}^{\si+i\infty}\tilde{f}_0(s)e^{-\be r \sqrt{s}}e^{st}\, \mathrm{d}s, \label{eq:calT1raph}\\
\mathcal{T}_2(r,\be,t) &=\frac{1}{2i\pi}\int_{\si-i\infty}^{\si+i\infty}\frac{\tilde{g}_0(s)}{\sqrt{s}}e^{-\be r \sqrt{s}}e^{st} \label{eq:calT2raph}
\end{align}
and where as usual $\si\in\mathbb{R}$ is chosen here such that all singularities of the integrands are to the left of the line $s=\si$.

\subsection{The case of $\theta\in(0,\pi/2]$} \label{sub:postheta}

This case follows entirely analogously to the negative $\theta$ scenario, the difference here being that the paths $B_{\pm}$ reside in the lower half of the complex $\al$-plane and as such the pole that leads to a contribution to the integral is at $\al=-i\eps$. We find that
%
\begin{align}
\cT(r,\theta,t) &= \frac{T_0}{\sqrt{2}\pi} \int_1^{\infty} \mathcal{G}(\be,\theta) \cT_1(r,\be,t) \, \mathrm{d}\be
+\frac{T_0'}{\sqrt{2}\pi}\int_1^{\infty} \mathcal{G}(\be,\theta) \cT_2(r,\be,t) \, \mathrm{d}\be \nonumber \\
& \hspace{4cm} + T_0 \cT_1(r,\cos\theta,t).
\end{align}

As an aside, for $\theta\in(0,\pi/2]$, analogously to the derivation of \eqref{eq:easynegraph} we have
\begin{align}
\int_1^{\infty}\mathcal{F}(\be,\theta)\, \mathrm{d}\be &= 0  \label{eq:intFzero}
\end{align}
Therefore using, \eqref{eq:easynegraph}, \eqref{eq:calGraph} and \eqref{eq:intFzero} we have
\begin{align}
\frac{1}{\pi\sqrt{2}}\int_1^{\infty}\mathcal{G}(\be,\theta)\, \mathrm{d}\be&=1- H(\theta)
\label{eq:magicidentity}
\end{align}
where $H(\theta)$ is the Heaviside step function.


\subsection{A summary of the solution}

Combining the results from sections \ref{sec:negtheta} and \ref{sub:postheta}, we can write down the solution for all values of $\theta$,
\begin{align}
\cT(r,\theta,t) 
&= \frac{T_0}{\pi \sqrt{2}} \int_1^{\infty} \mathcal{G} (\beta, \theta)
  \mathcal{T}_1 (r, \beta, t) \mathrm{d} \beta + \frac{T_0'}{\pi \sqrt{2}}
  \int_1^{\infty} \mathcal{G} (\beta, \theta) \mathcal{T}_2 (r, \beta, t)
  \mathrm{d} \beta \nonumber\\
  & + H(\theta) T_0 \mathcal{T}_1 (r, \cos \theta, t) - (1
  - H(\theta)) T_0' \mathcal{T}_2 (r, \cos \theta, t).
\label{eq:generalsolution1}
\end{align}
Both integrals in \eqref{eq:generalsolution1} and the additional term have a discontinuity at $\theta=0$ and as such this form is not particularly ``clean''. We are able to improve upon this form, using \eqref{eq:magicidentity} to generate
\begin{multline}
\cT(r,\theta,t) = \frac{T_0}{\pi \sqrt{2}} \int_1^{\infty} \mathcal{G} (\beta, \theta)
  \{ \mathcal{T}_1 (r, \beta, t) -\mathcal{T}_1 (r, \cos \theta, t) \}
  \mathrm{d} \beta \\
  + \frac{T_0'}{\pi \sqrt{2}} \int_1^{\infty} \mathcal{G} (\beta, \theta)
  \{ \mathcal{T}_2 (r, \beta, t) -\mathcal{T}_2 (r, \cos \theta, t) \}
  \mathrm{d} \beta
  + T_0 \mathcal{T}_1 (r, \cos \theta, t).  \label{eq:niceexpre}
\end{multline}
Each term in this expression is now continuous across $\theta=0$. We shall discuss this aspect further in the context of specific examples in the next section.

\section{Some specific boundary conditions} \label{sec:particularexamples}

\subsection{Perfect insulator on $y<0$}

Let us now assume that $T_0'=0$ so that we have a perfect insulator on $y<0$. We shall consider a variety of temperature profiles for $y>0$. The solution is therefore obtained by setting $T_0'=0$ in \eqref{eq:generalsolution1} or equivalently \eqref{eq:niceexpre}. As such only $\cT_1$ enters the analysis.


\subsubsection{Step temperature change}\label{sub:stepchange}

\begin{figure}[h!]
\psfrag{y}{$y$}
\psfrag{T}{$\cT$}
\begin{center}
\includegraphics[scale=0.9]{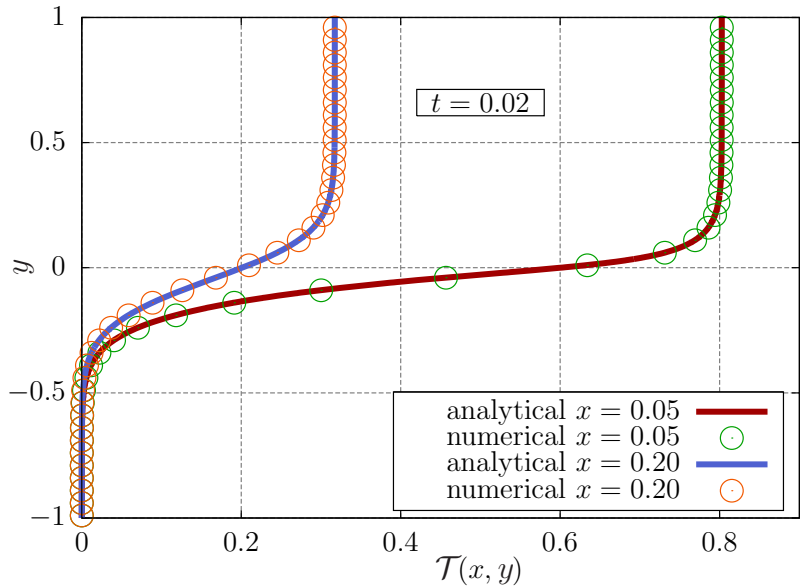}
\end{center}
\caption{Thermal field $\cT(x,y,t)$ at $x=0.05$ (red) and $x=0.2$ (blue) when $t=0.02$ for the case of a perfect insulator on $y<0$ and a step change temperature increase at $t=0$ on $y>0$. Circles are predictions from solutions to the same problem using finite elements methods in COMSOL, which provides validation of the present semi-analytical method.}
\label{localcoords}
\end{figure}

Take the simplest form, $f_0(t)=1$ so that $\tilde{f}_0(s)=1/s$ and we need to determine $\cT_1$ defined in \eqref{eq:calT1raph}. It transpires that it is convenient to differentiate the expression for $\cT_1$ with respect to $t$, which enables the inverse Laplace integral to be evaluated analytically in this case
\begin{align}
\frac{\mathrm{d}\cT_1}{\mathrm{d}t}(r,\be,t) &= \frac{1}{2\pi i}\int_{\si-i\infty}^{\si+i\infty}e^{-\be r\sqrt{s}}e^{st}\lit ds = \frac{r\be}{2\sqrt{\pi}t^{3/2}}e^{-r^2\be^2/(4t)}.
\end{align}
We then integrate (definitely) with respect to $t$ with a lower limit of $t=0$, finding that
\begin{align}
\cT_1(r,\be,t) &= \tn{Erfc}\left(\frac{r\be}{2\sqrt{t}}\right).
\label{eq:Tcal1_simple}
\end{align}
Appealing to \eqref{eq:generalsolution1} we have as our solution
\begin{align}
\cT(r,\theta,t) &= \frac{T_0}{\sqrt{2}\pi}\int_1^{\infty} \mathcal{G}(\be,\theta) \tn{Erfc}\left(\frac{r\be}{2\sqrt{t}}\right) \, \mathrm{d}\be
+H(\theta)T_0 \tn{Erfc}\left(\frac{r\cos\theta}{2\sqrt{t}}\right)
\label{eq:discontinuous_step_T}
\end{align}
noting that $G(\be,\theta)$ is the real function defined in \eqref{calG}. Each term of \eqref{eq:discontinuous_step_T} is easily computed numerically, and, in Fig.\ \ref{localcoords} we plot the resulting temperature profile on the horizontal axis against $y$ running vertically, at time $t=0.02$ for two values of $x$, $x=0.05$ and $x=0.2$. The circles are results taken from a finite element solution of the same problem in COMSOL and provide validation of the present semi-analytical scheme. 

\begin{figure}[h!]
\centerline{
\includegraphics[width=0.5\textwidth]{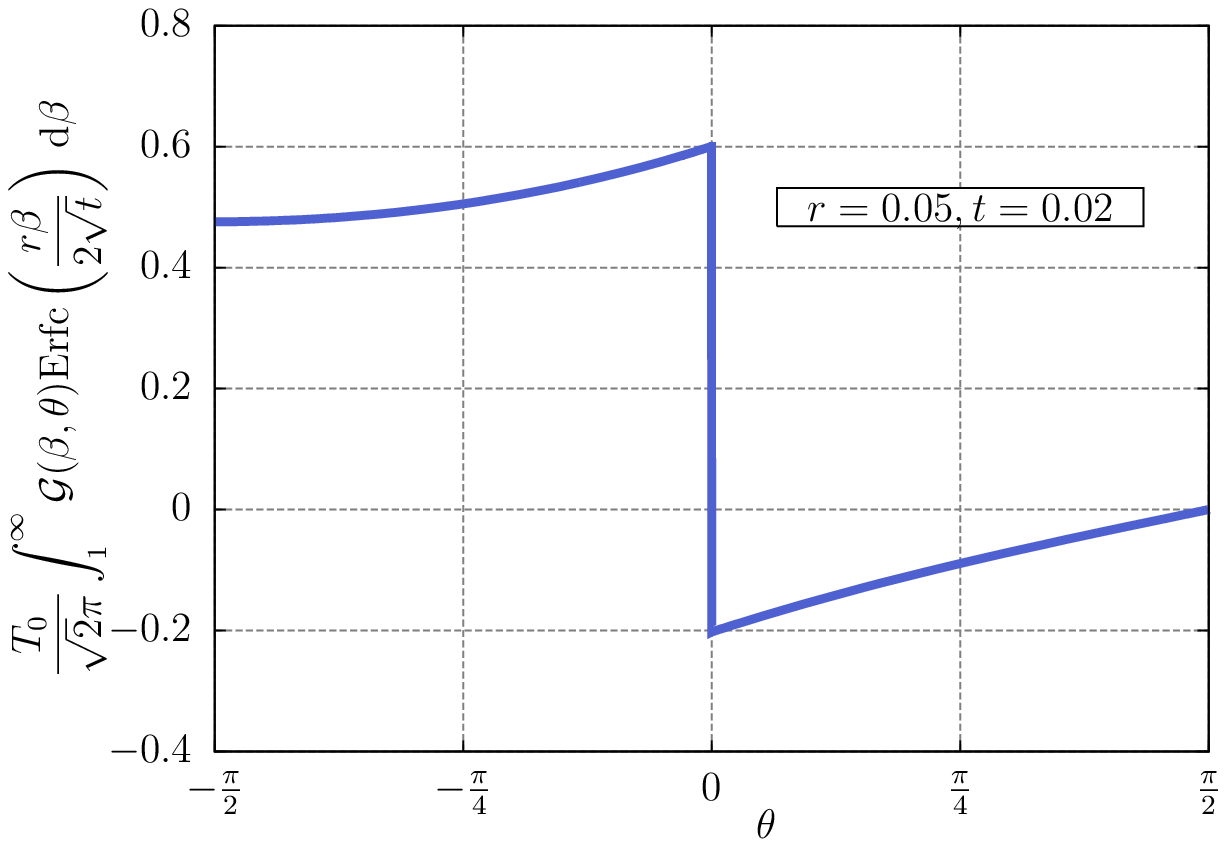}
\includegraphics[width=0.5\textwidth]{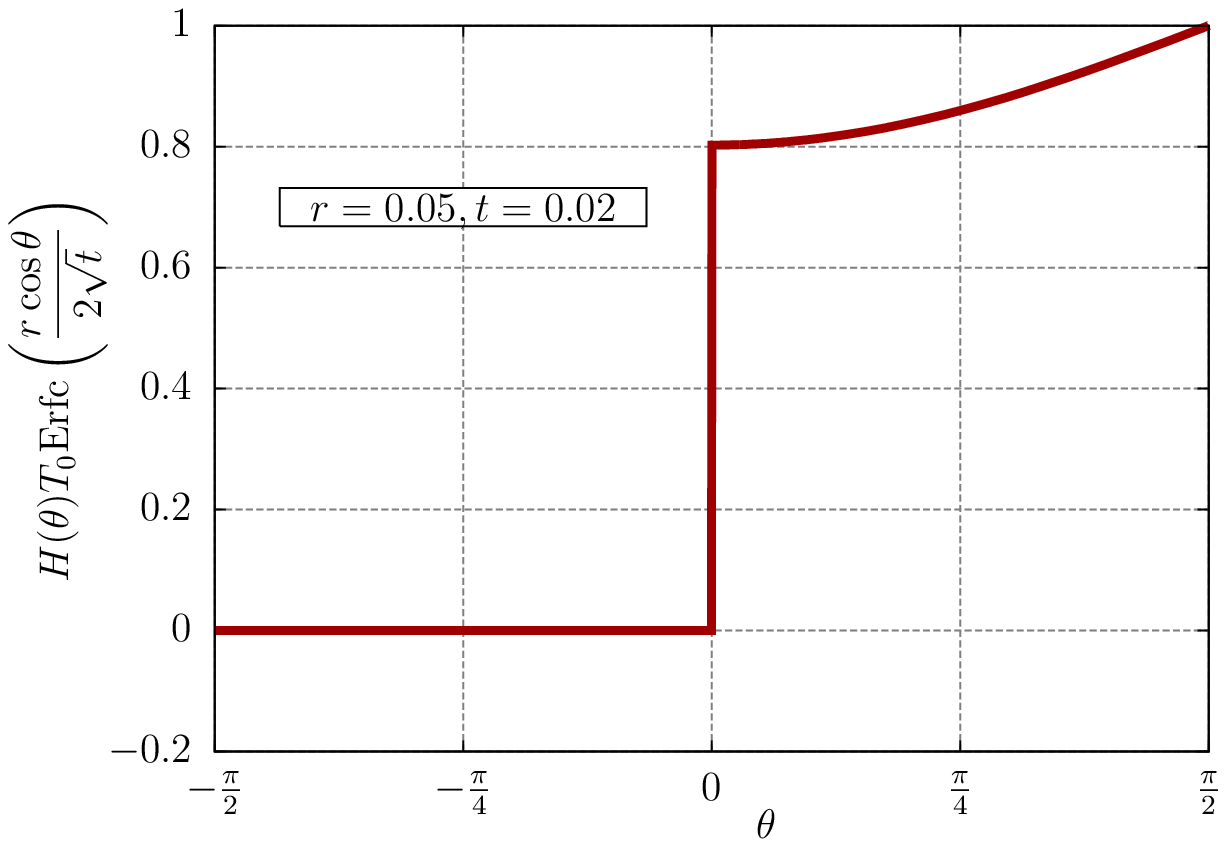}}
\caption{Plot of the first (left) and
  second (right) terms on the right hand side of \eqref{eq:discontinuous_step_T} for $r=0.05$, $t=0.02$ and $\theta \in [- \pi / 2,
  \pi/2]$.}
\label{fig:discontinuous_terms}
\end{figure}

\begin{figure}[h!]
 \centerline{\includegraphics[width=0.5\textwidth]{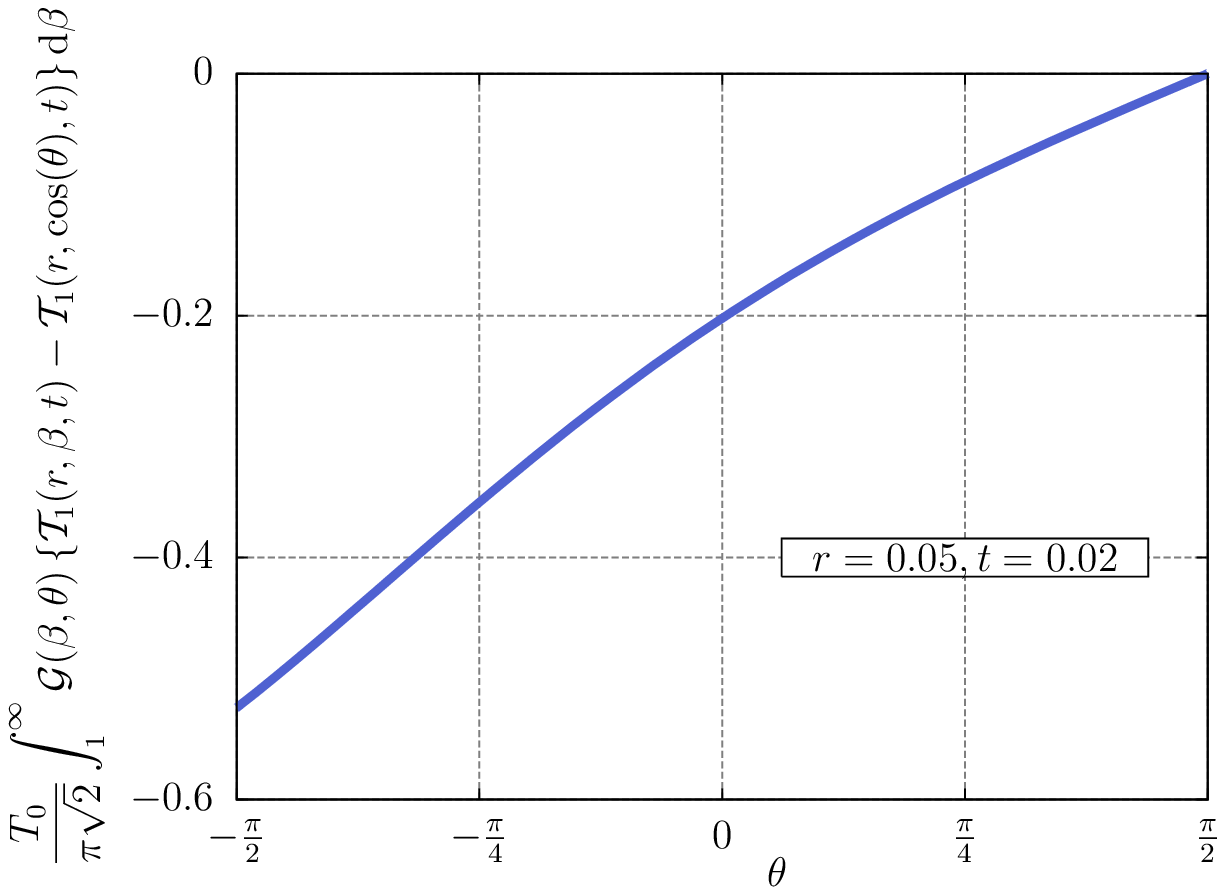}\includegraphics[width=0.5\textwidth]{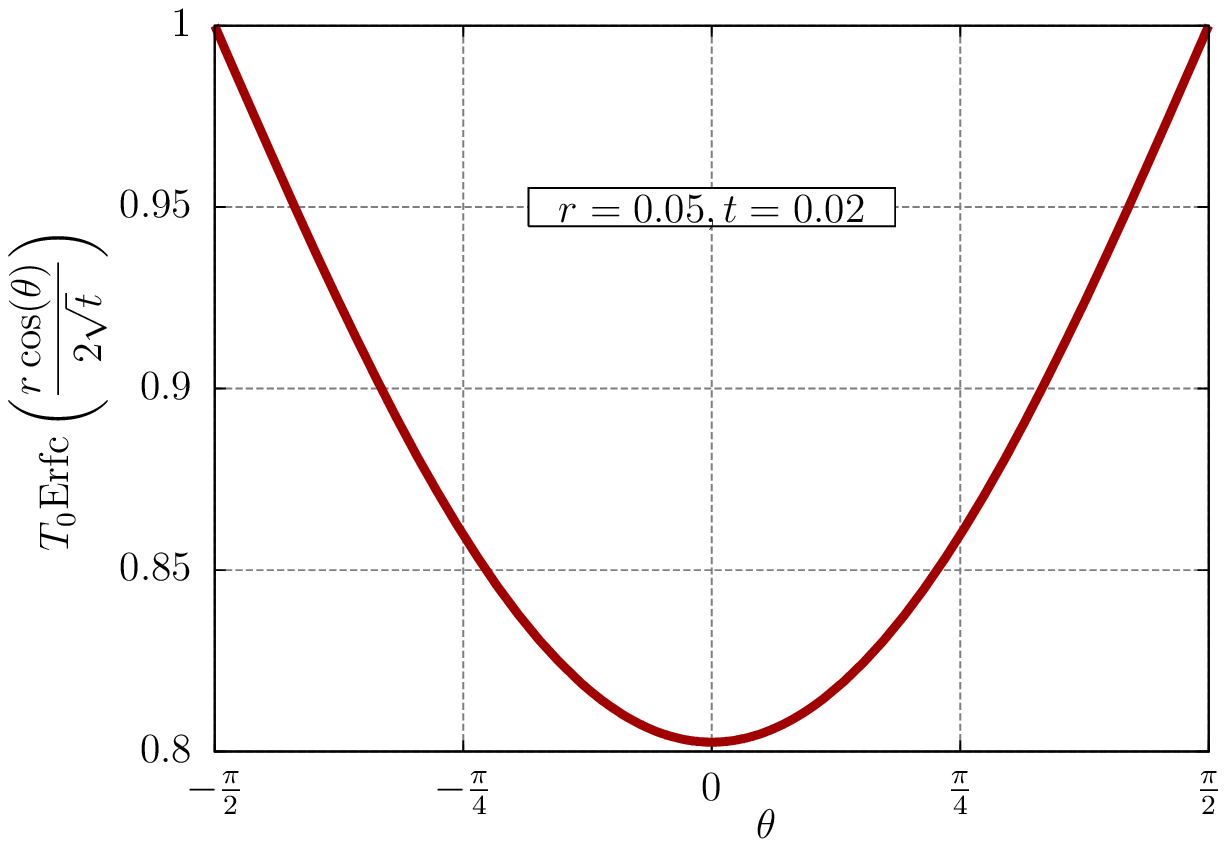}}
  \caption{Plot of the first (left) and
  second (right) terms on the right hand side of \eqref{eq:continuous_step_T} for $r=0.05$, $t=0.02$ and $\theta \in [- \pi / 2,
  \pi/2]$.}
\label{fig:continuous_terms}
\end{figure}

Note that the temperature profile is continuous across the $x$-axis, i.e.\ $\theta=0$. Rather interestingly, both terms on the right hand side of \eqref{eq:discontinuous_step_T} are discontinuous across $\theta=0$, as is shown in Fig.\ \ref{fig:discontinuous_terms} but the two discontinuities compensate exactly to yield a continuous temperature profile.
If instead of the form \eqref{eq:generalsolution1}, we use \eqref{eq:niceexpre}, the solution is written as
\begin{multline}
\cT(r,\theta,t) =\frac{T_0}{\pi \sqrt{2}} \int_1^{\infty} \mathcal{G} (\beta, \theta)
  \left\{ \tn{Erfc}\left(\frac{r\be}{2\sqrt{t}}\right) -\tn{Erfc}\left(\frac{r\cos(\theta)}{2\sqrt{t}}\right)\right\}
  \mathrm{d} \beta\\
  +T_0 \tn{Erfc}\left(\frac{r\cos(\theta)}{2\sqrt{t}}\right).
\label{eq:continuous_step_T}
\end{multline}
This expression is also straightforward to evaluate numerically and (obviously) gives the same results as those presented in Fig.\ \ref{localcoords}, but this time, as shown in Fig.\ \ref{fig:continuous_terms}, both terms on the right hand side of \eqref{eq:continuous_step_T} are continuous across $\theta=0$.

\begin{figure*}[h!]\centering
\includegraphics[scale=1.0]{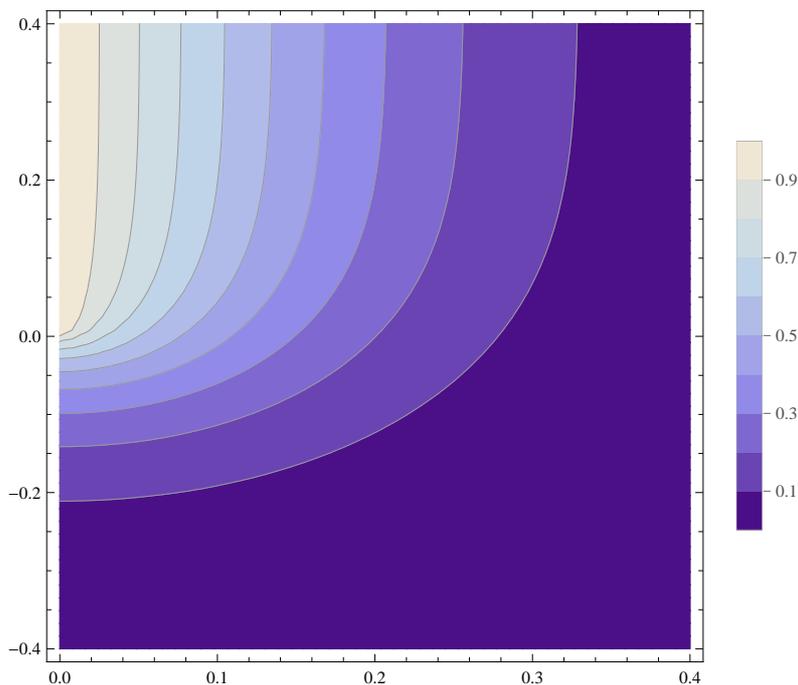}
\caption{Contour plot of the temperature profile when $t=0.02, T_0=1, T_0'=0$, i.e.\ a perfect insulator on $y<0$ and step change in temperature on $y>0$. We see how the thermal field propagates into $y<0$.}
\label{fig:temp1}
\end{figure*}

Finally, in Fig.\ \ref{fig:temp1}, we plot the two dimensional temperature contour profile on the $(x,y)$ plane at $t=0.02$ illustrating how the distribution spreads out from the upper half-plane.

\subsection{Imperfect insulator on $y<0$}

\begin{figure}[h!]
\psfrag{y}{$y$}
\psfrag{T}{$\cT$}
\begin{center}
\includegraphics[scale=0.9]{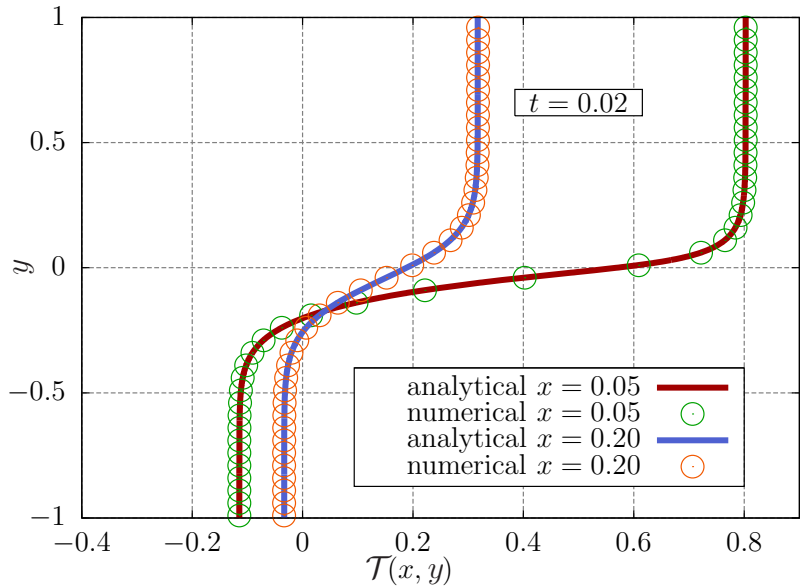}
\end{center}
\caption{Thermal field $\cT(x,y,t)$ at $x=0.05$ (red) and $x=0.2$ (blue) when $t=0.02$ for the case of an imperfect insulator on $y<0$ and a step change temperature increase at $t=0$ on $y>0$. Circles are predictions from solutions to the same problem using finite elements methods in COMSOL, which provides validation of the present semi-analytical method.}
\label{fig:tempprof5.2}
\end{figure}

\begin{figure*}[hbt!]\centering
\includegraphics{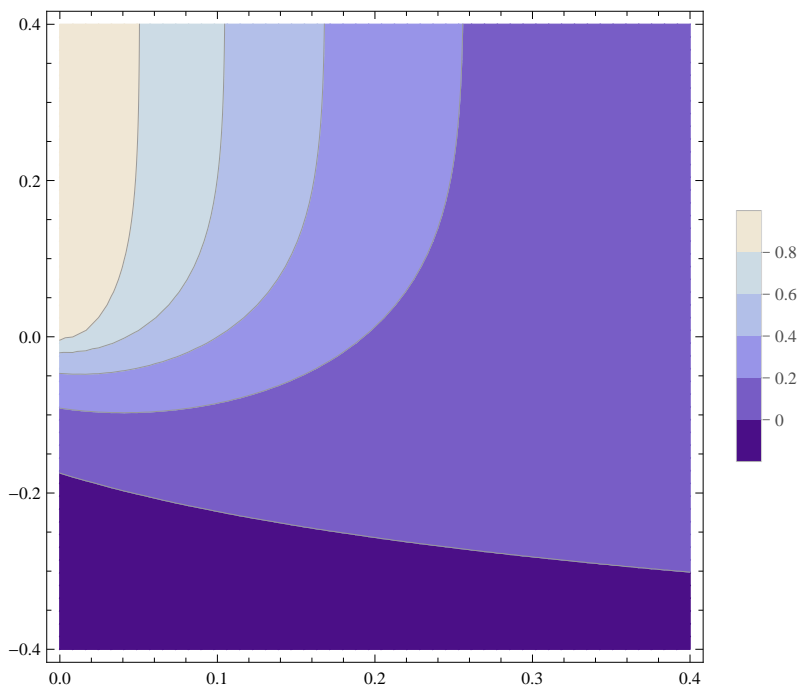}
\caption{Contour plot of the temperature profile when $t=0.02, T_0=1, T_0'=1$, i.e.\ an imperfect insulator on $y<0$ and step change in temperature on $y>0$. We see how the thermal field propagates into $y<0$.}
\label{fig:temp2}
\end{figure*}

Let us now consider the case when $T_0'\neq 0$ and let us take $f(t)=g(t)=1$ so that $\tilde{f}_0(s)=\tilde{g}_0(s)=1/s$ and now both $\cT_1$ and $\cT_2$ play a role. Once again it is convenient to differentiate with respect to $t$ in order to evaluate the inverse Laplace transforms
and subsequently integrating these expressions definitely with respect to $t$ with a lower limit of $t=0$ yields
\begin{align}
\cT_1(r,\be,t) &= \tn{Erfc}\left(\frac{r\be}{2\sqrt{t}}\right), \\
\cT_2(r,\be,t) &= 2\sqrt{\frac{t}{\pi}} e^{-\be^2r^2/(4t)}- r \be \tn{Erfc}\left(\frac{r\be}{2\sqrt{t}}\right).
\end{align}
Either of the expressions \eqref{eq:generalsolution1} or \eqref{eq:niceexpre} then recover the temperature profile. Both formulations are easily computed numerically and give rise to a continuous temperature profile, as seen in Fig.\ \ref{fig:tempprof5.2} where the profile is plotted at $t=0.02$ for $x=0.05$ and $x=0.2$. The contour plot of the thermal field at $t=0.02$ is given in Fig.\ \ref{fig:temp2}.

\subsection{Continuous ramp up and down superposed on a step change in $y>0$}

\begin{figure}[h]
\psfrag{f}{$f_0(t)$}
\psfrag{t}{$t$}
\psfrag{0}{$0$}
\psfrag{a}{$a$}
\psfrag{b}{$b$}
\psfrag{d}{$2b-a$}
\psfrag{1}{$1$}
\psfrag{c}{$1+b-a$}
\psfrag{e}{$0$}
  \centering\includegraphics[width=0.8\textwidth]{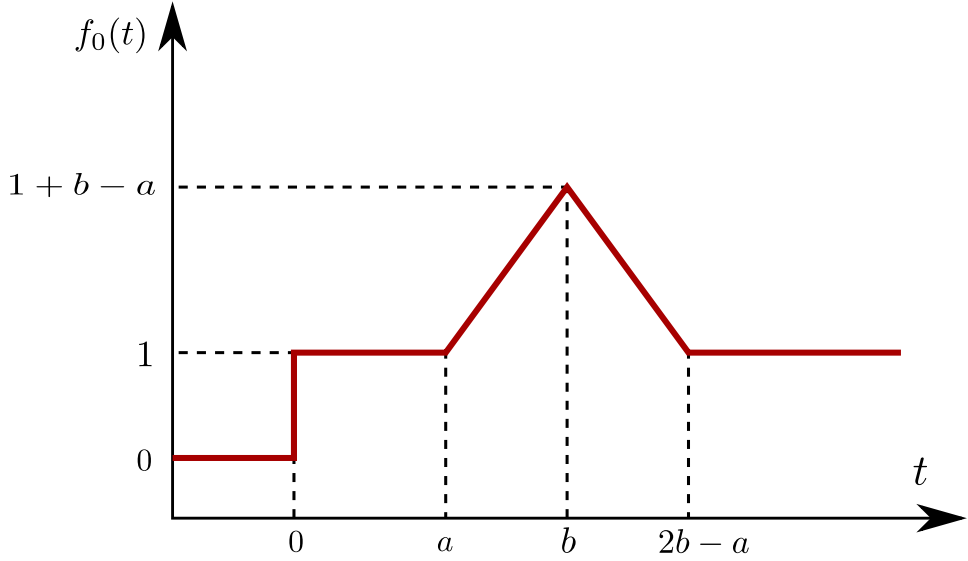}
  \caption{Plot of $f_0 (t)$ as defined by \eqref{f0ramp}. A ramp up and down is superposed on a step change unit temperature profile.}
	\label{fig:rampupdown}
\end{figure}

Thus far we have considered only cases where $f_0$ and $g_0$ are constant, accommodating for the step change at $t=0$ of course. Let us now consider the case when these functions can be unsteady and in particular when $f_0(t)$ is a step and superposed general ramp up and down profile given by
\begin{align}
  f_0 (t) & =  H(t) + (t - a) H(t - a) + 2 (b - t)
  H(t - b) + (t - (2 b - a)) H(t - (2 b - a)) \label{f0ramp}
\end{align}
and illustrated in Fig.\ \ref{fig:rampupdown}.


%

\begin{figure}[h!]
  \centerline{\includegraphics[width=0.5\textwidth]{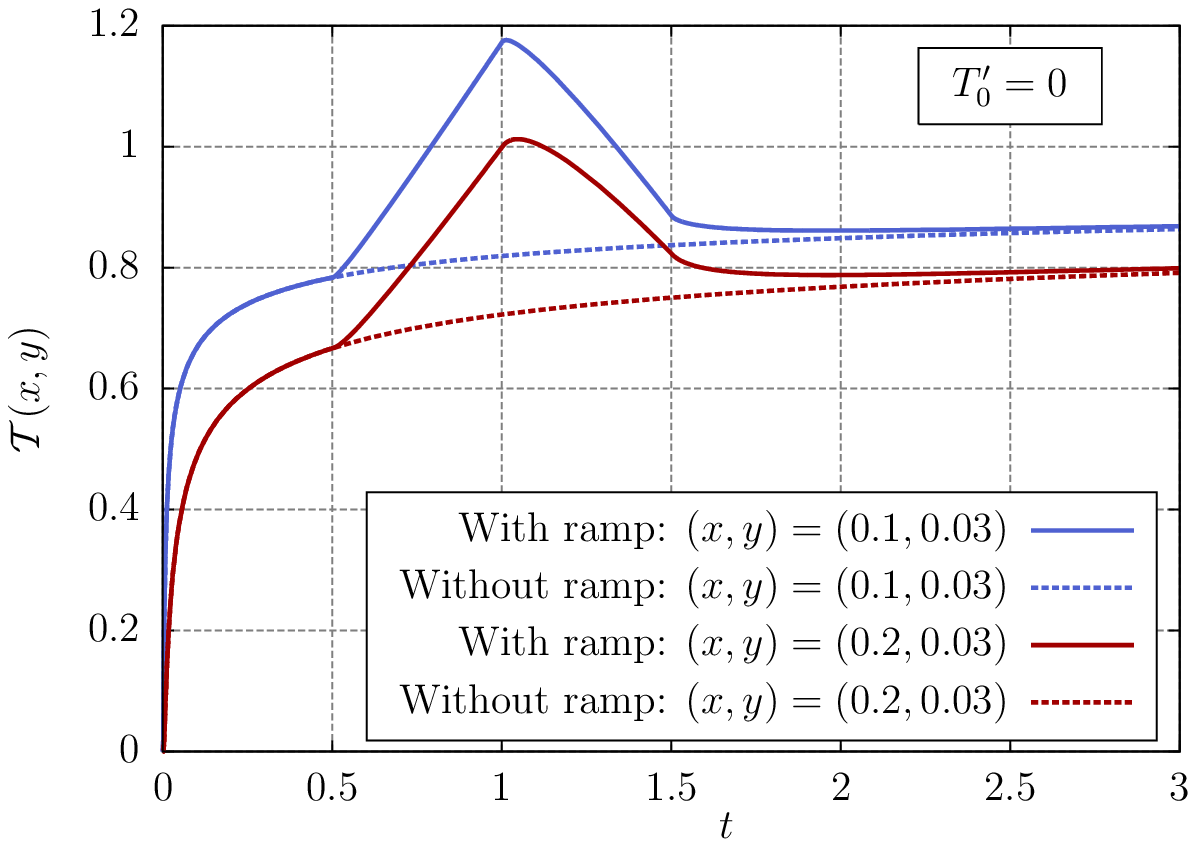}\quad\includegraphics[width=0.5\textwidth]{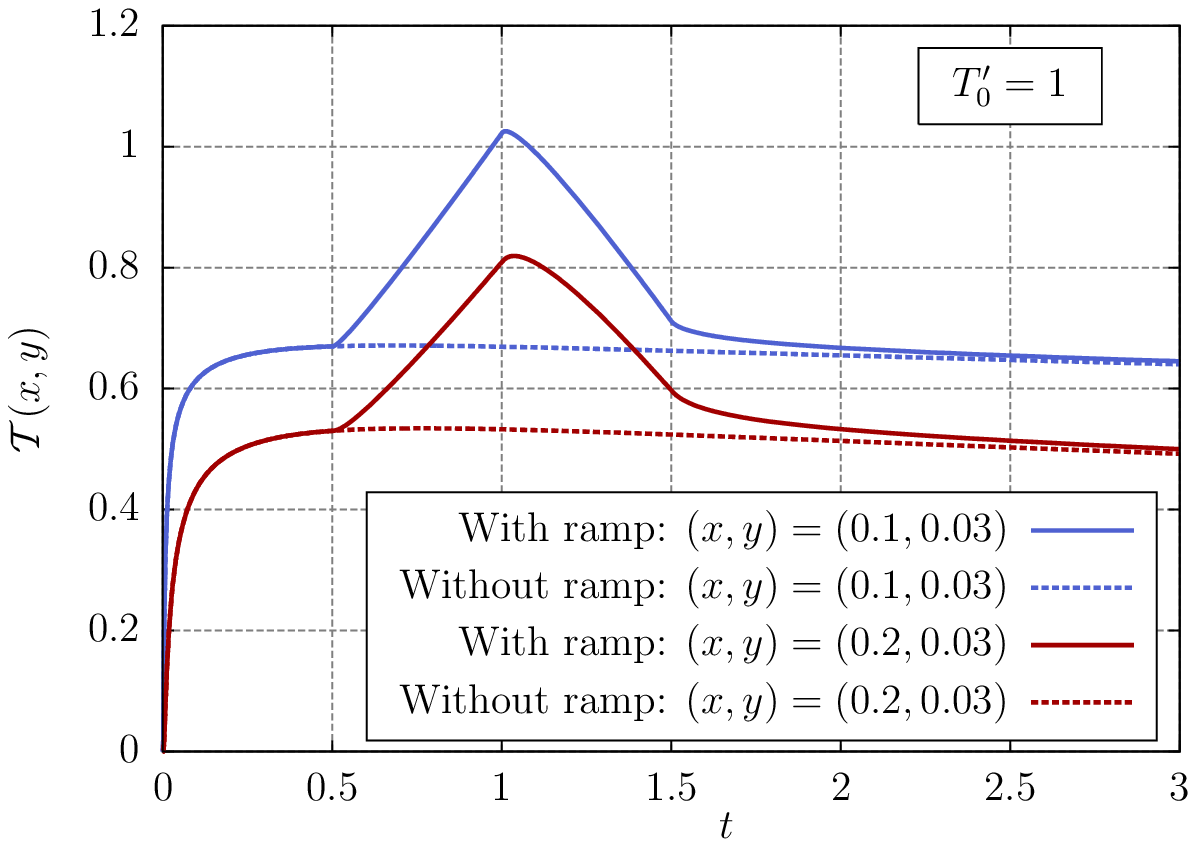}}
  \caption{Illustrating the effect of the ramp up and down on the evolution of the temperature profile at
  two different locations $(x, y) = (0.1, 0.03)$ (red) and $(x, y) = (0.2,
  0.03)$ (blue) for $T_0 = 1$ and $T_0'=0$ (left) or $T_0 = 1$ and $T_0' = 1$ (right). The
  dashed lines correspond to the profile without the unsteady ramp, and the
  plain lines correspond to the profile with the unsteady ramp.}
  \label{rampfigs}
\end{figure}

The Laplace transform of $f_0 (t)$ is
\begin{align}
  \tilde{f}_0 (s) & =  \frac{1}{s} + \frac{e^{- as}}{s^2} - 2 \frac{e^{-
  bs}}{s^2} + \frac{e^{- (2 b - a) s}}{s^2} \nonumber \\
  & =  \tilde{f}_0^{(1)} (s) + \tilde{f}_0^{(2)} (s) + \tilde{f}_0^{(3)} (s)
  + \tilde{f}_0^{(4)} (s)
\end{align}
and the resulting expression for $\mathcal{T}_1$ is
\begin{align}
  \mathcal{T}_1 (r, \beta, t) & = \frac{1}{2 i \pi} \int_{\si - i \infty}^{\si +
  i \infty} (\tilde{f}_0^{(1)} (s) + \tilde{f}_0^{(2)} (s) + \tilde{f}_0^{(3)}
  (s) + \tilde{f}_0^{(4)} (s)) e^{- \beta r \sqrt{s}} e^{st} \mathd s \nonumber \\
  & =  \mathcal{T}_1^{(1)} (r, \beta, t) +\mathcal{T}_1^{(2)} (r, \beta, t)
  +\mathcal{T}_1^{(3)} (r, \beta, t) +\mathcal{T}_1^{(4)} (r, \beta, t).
\end{align}
The case of $\mathcal{T}_1^{(1)}$ has already been dealt with in Section \ref{sub:stepchange}
and is thus given by \eqref{eq:Tcal1_simple}. Moreover, $\mathcal{T}_1^{(2)} (r,
\beta, t), \hspace{1em} \mathcal{T}_1^{(3)} (r, \beta, t)$ and
$\mathcal{T}_1^{(4)} (r, \beta, t)$ have a very similar structure, so we only
need to consider only the case of $\mathcal{T}_1^{(2)} (r, \beta, t)$ in detail. In fact, we have
\begin{align}
  \mathcal{T}_1^{(2)} (r, \beta, t) & =  \frac{1}{2 i \pi}\int_{\si - i \infty}^{\si + i \infty}
  \tilde{f}_0^{(2)} (s) e^{- \beta r \sqrt{s}} e^{st} \mathd s \nonumber \\
  & =   \frac{1}{2 i \pi}\int_{\si - i \infty}^{\si + i \infty} \frac{e^{- as}}{s^2} e^{- \beta r
  \sqrt{s}} e^{st} \mathd s
\end{align}
Differentiating $\mathcal{T}_1^{(2)}$ twice with respect to time, we obtain
\begin{align}
  \frac{\partial^2 \mathcal{T}_1^{(2)}}{\partial t^2} (r, \beta, t) & =
  H(t - a) \frac{r \beta}{2 \sqrt{\pi} (t - a)^{3 / 2}} e^{-
  \frac{r^2 \beta^2}{4 (t - a)}},
\end{align}
which, following the same reasoning as in Section \ref{sub:stepchange} implies that
\begin{align}
  \frac{\partial \mathcal{T}_1^{(2)}}{\partial t^{}} (r, \beta, t) & =
  H(t - a) \tmop{Erfc} \left( \frac{r \beta}{2 \sqrt{t - a}}
  \right)
\end{align}
and by definite integration with respect to time, we obtain
\begin{align}
  \mathcal{T}_1^{(2)} (r, \beta, t) & =  H(t - a) \left\{ \left( t
  - a + \frac{r^2 \beta^2}{2} \right) \tmop{Erfc} \left( \frac{r \beta}{2
  \sqrt{t - a}} \right) - \frac{r \beta \sqrt{t - a}}{\sqrt{\pi}} e^{-
  \frac{r^2 \beta^2}{4 t}} \right\}.
\end{align}
The functions $\mathcal{T}_1^{(3)} (r, \beta, t)$ and $\mathcal{T}_1^{(4)} (r,
\beta, t)$ are obtained in the same manner. Using these expressions in the general formulation \eqref{eq:generalsolution1} or \eqref{eq:niceexpre} enables the solution to be computed rather rapidly. In Fig.\ \ref{rampfigs} we plot the resulting thermal field at the locations $(x,y)=(0.1,0.03)$ and $(0.2,0.03)$ in the cases when $T_0=1$ with $T_0'=0$ (left) and $T_0=1$ with $T_0'=1$ (right). We also plot the associated solutions determined previously where no ramp up and down is present in the boundary condition.

%
%
%
%

\section{Concluding remarks} \label{sec:conc}

Employing the Wiener-Hopf technique and a Cagniard-de Hoop-type integral, a rapidly convergent integral expression has been determined for a class of transient thermal mixed boundary value problems. The integral is easily computable on a standard desktop PC for a wide range of transient boundary forcings of interest. Here we illustrated the computation for a number of cases, including step-changes in temperature and ramp up and down boundary profiles. Such quasi-analytical expressions are of great utility in order to speed up computations and enable asymptotic analysis close to locations of interest. Future work could include extensions to full elastodynamics and coupled thermoelasticity. In these cases matrix Wiener-Hopf problems will result in general.

\subsection*{Acknowledgements}

Parnell is thankful to Universit\'{e} Paris-Est Cr\'eteil (UPEC) for funding his visiting position in April 2013 when this work was initiated. He also gratefully acknowledges the Engineering and Physical Sciences Research Council for funding his research fellowship (EP/L018039/1). Abrahams thanks the Royal Society for a Wolfson Research Merit award (2013-2018).


\appendix

\section{Simplification of the function $\mathcal{F}(\beta,\theta)$ } \label{appendix}

\begin{figure}[h!]
\psfrag{Bm}{$B_-$}
\psfrag{Bp}{$B_+$}
\psfrag{pm}{$\psi_-$}
\psfrag{pp}{$\psi_+$}
\psfrag{rp}{$r_+=R$}
\psfrag{rm}{$r_-=R$}
\begin{center}
\includegraphics[scale=0.6]{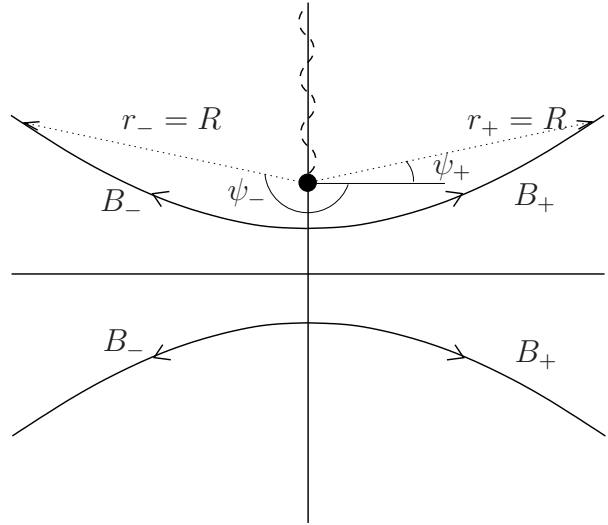}
\end{center}
\caption{Diagrammatic description of the angles $\psi^\pm$}
\label{fig:localcoords}
\end{figure}

From \eqref{eq:calFraph} we have
\begin{align}
\mathcal{F}(\be,\theta) &= \frac{1}{\al_+(\al_+-i)^{1/2}}\frac{d\al_+}{d\be} -
\frac{1}{\al_-(\al_- -i)^{1/2}}\frac{d\al_-}{d\be}
\end{align}
and
\begin{align}
\frac{d\al_{\pm}}{d\be} &= -i\sin\theta \pm \frac{\be}{\sqrt{\be^2-1}}\cos\theta.
\end{align}
Therefore
\begin{multline}
\mathcal{F}(\be,\theta) = \frac{1}{\al_+(\al_+-i)^{1/2}}\left(-i\sin\theta + \frac{\be}{\sqrt{\be^2-1}}\cos\theta\right)\\
+
\frac{1}{\al_-(\al_- -i)^{1/2}}\left(i\sin\theta + \frac{\be}{\sqrt{\be^2-1}}\cos\theta\right)
\end{multline}
and simplifying further we obtain
\begin{multline}
\mathcal{F}(\be,\theta) = \frac{1}{\al_-\al_+(\al_+-i)^{1/2}(\al_--i)^{1/2}}\Bigg[
(\al_+-i)^{1/2}\left(\be-\frac{i\cos\theta\sin\theta}{\sqrt{\be^2-1}}\right)\\
-
(\al_--i)^{1/2}\left(\be+\frac{i\cos\theta\sin\theta}{\sqrt{\be^2-1}}
\right)\Bigg].
\end{multline}
Next with reference to Fig.\ \ref{fig:localcoords} define $R$ and $\psi_{\pm}$ such that
\begin{align}
\al_+-i &= R e^{i\psi_+}, & \al_--i &= R e^{i\psi_-}
\end{align}
and since $\psi_-=-\pi-\psi_+$ we can write
\begin{align}
(\al_+-i)^{1/2} &= \sqrt{R}e^{i\psi_+/2}, & (\al_--i)^{1/2} &= \sqrt{R}e^{-\psi_-/2} = -i\sqrt{R}e^{-\psi_+/2}.
\end{align}
Therefore it is possible to show that
\begin{multline}
\mathcal{F}(\be,\theta) 
=  \frac{1}{\al_-\al_+(\al_+-i)^{1/2}(\al_--i)^{1/2}}\sqrt{R}(1+i)
\Big[\be\left(\cos\left(\psi_+/2\right)+\sin\left(\psi_+/2\right)\right)-\\
\cos\theta\sin\theta(\be^2-1)^{-1/2}\left(\cos\left(\psi_+/2\right)-\sin\left(\psi_+/2\right)\right)\Big]. \label{FFF}
\end{multline}
Finally, noting that
\begin{align}
\al_+\al_- = -(\be^2-\cos^2\theta)
\end{align}
and
\begin{align}
(\al_+-i)^{1/2}(\al_--i)^{1/2} 
&= -i|\be+\sin\theta|
\end{align}
so that
\begin{align}
\al_-\al_+(\al_+-i)^{1/2}(\al_--i)^{1/2} = i (\be^2-\cos^2\theta)|\be+\sin\theta|
\end{align}
we have from \eqref{FFF}
\begin{align}
\mathcal{F}(\be,\theta) &= \frac{\sqrt{2}}{ic} \mathcal{G}(\be,\theta)
\end{align}
where
\begin{multline}
\mathcal{G}(\be,\theta) = \frac{1}{{(\be^2-\cos^2\theta)|\be+\sin\theta|^{1/2}}}
\Big[\be\left(\cos\left(\psi_+/2\right)+\sin\left(\psi_+/2\right)\right)-\\
\cos\theta\sin\theta(\be^2-1)^{-1/2}\left(\cos\left(\psi_+/2\right)-\sin\left(\psi_+/2\right)\right)\Big] \label{calG}
\end{multline}
is a real valued function.

\vspace{1cm}
\normalsize

\section{Variational formulation used for finite element solution}

In order to compute the transient solution $\mathcal{T}(x,y,t)$ of the problem \eqref{2.4}-\eqref{2.6} in the semi-infinite domain $\mathcal{D}$ using finite element method, we define a rectangular domain $\mathcal{D}_0 \in \mathbb{R}^2$: $\mathcal{D}_0 = [0,a]\times [-b/2,b/2]$ with boundary $\partial{D}_0 = \partial{D}_0^-\cup\partial{D}_0^+\cup\partial{D}_0^\infty$ where $\partial{D}_0^- = \{x=0,-b/2<y<0\}$ and $\partial{D}_0^+ = \{x=0,0<y<b/2\}$. The parameters $a$, $b$ which define the size of $\mathcal{D}_0$ are chosen sufficiently large such that the temperature field on $\partial{D}_0^\infty$ is not influenced by the perturbation due to the discontinuous boundary condition on $\partial{D}_0^+$ and $\partial{D}_0^-$ in the time interval under consideration.
As a consequence, the symmetric boundary condition may be imposed on $\partial{D}_0^\infty$. The problem \eqref{2.4}-\eqref{2.6}
may be rewritten as:
\begin{align}
    & \frac{\partial\mathcal{T}}{\partial t}
       -\Grad^2{\mathcal{T}} = 0, & \text{in}\; \partial{D}_0\; \\
    & \mathcal{T} = T_0f_0(t),  &
    	\text{on}\;  \partial\mathcal{D}_0^+ \\
    & \Grad{\mathcal{T}}\cdot\vc{n} =
    		-T'_0g_0(t), & \text{on}\; \partial{D}_0^- \\
	& \Grad{\mathcal{T}}\cdot\vc{n} = 0, & \text{on}\;
		\partial{D}_0^\infty
\end{align}
where $\vc{n}$ denotes the outward pointing normal vector to $\partial{D}_0$.

Let $\dT$ be the trial function of $\mathcal{T}$,
the weak formulation of this problem reads: For any test function $\dT \in \bar{V}(\mathcal{D}_0)$, find $\mathcal{T} \in V(\mathcal{D}_0)$ such that
\begin{align}
  & \int_{\mathcal{D}_0} \dT~\frac{\partial\mathcal{T}}{\partial t}\,dv +
  \int_{\mathcal{D}_0} \Grad(\dT)\cdot\Grad{T}\,dv
  = -\int_{\partial\mathcal{D}_0^-} (\dT)\,T'_0g_0\,ds, \\
   & \text{and} \quad \mathcal{T} = T_0f_0 \quad \text{on}\; \partial\mathcal{D}_0^+,
\end{align}
where $V(\mathcal{D}_0) = \{f(x,y) \in H^1(\mathcal{D}_0)\}$,
$\bar{V}(\mathcal{D}_0) = \{f(x,y) \in H^1(\mathcal{D}_0); f(x,y) = 0, (x,y) \in \partial\mathcal{D}_0^+\}$.

This weak formulation has been implemented in the finite element software COMSOL Multiphysics \cite{Comsol2008}. For simulations
in the time interval  $t = [0,0.02]$ as presented in Section 5, a rectangular domain $(x,y) \in [0,1]\times[-1,1]$ has been used. The element size ($h_e$) and the time step ($\Delta t$) needed for the discretization are respectively
$h_e = 0.01$ and $\Delta{t} = 10^{-4}$.

\bibliography{thermo_refs}
\bibliographystyle{unsrt}

\end{document}

%% file: macros.tex
\newcommand{\al}{\alpha}

\newcommand{\be}{\beta}

\newcommand{\eps}{\epsilon}

\newcommand{\ka}{\kappa}

\newcommand{\si}{\sigma}

\newcommand{\pa}{\partial}

\newcommand{\been}{\begin{equation}}
\newcommand{\een}{\end{equation}}
\newcommand{\beena}{\begin{eqnarray}}
\newcommand{\eena}{\end{eqnarray}}

\newcommand{\tn}{\textnormal}

\newcommand{\deriv}[2]{\frac{\partial{#1}}{\partial{#2}}}
\newcommand{\derivtwo}[2]{\frac{\partial^{2}{#1}}{\partial{#2}^2}}

\newcommand{\lit}{\hspace{0.2cm}}

\newcommand{\cT}{\mathcal{T}}
\newcommand{\mathd}{\mathrm{d}}
\newcommand{\tmop}[1]{\ensuremath{\operatorname{#1}}}
\newcommand{\Grad}{\nabla}
\newcommand{\vc}[1]{\boldsymbol{#1}}
\newcommand{\dT}{\delta{\mathcal{T}}}